\input amstex.tex
\documentstyle{amsppt}
\magnification=1100
\hsize = 6.0 truein
\vsize = 9.0 truein
\pageno 1
\topmatter
\title On the same $N$-type of the suspension of the infinite quaternionic projective space \endtitle
\leftheadtext{Dae-Woong Lee} \rightheadtext{On the same $N$-type of the suspension}
\author Dae-Woong Lee
\endauthor
\thanks
\newline
This paper was supported by research funds of Chonbuk National University in 2009.
\endthanks

\abstract
Let $[\rho_{i_k},[\rho_{i_{k-1}},\ldots,[\rho_{i_{1}}, \rho_{i_2}] \ldots]]$ be an iterated commutator of self-maps $\rho_{i_j} : \Sigma {\Bbb H}P^\infty \rightarrow \Sigma {\Bbb H}P^\infty , j = 1,2, \ldots, k$ on the suspension of the infinite quaternionic projective space.
In this paper, it is shown that the image of the homomorphism induced by the adjoint of this commutator is both primitive and decomposable.
The main result in this paper asserts that the set of all homotopy types of spaces having the same $n$-type as the suspension of the infinite quaternionic projective space is the one element set consisting of a single homotopy type.
Moreover, it is also shown that the group $\text{Aut}(\pi_{\leq n} (\Sigma {\Bbb H}P^\infty )/\text{torsion})$ of automorphisms is finite for $n \leq 9$, and infinite for $n \geq 13$, and that $\text{Aut}(\pi_{*} (\Sigma {\Bbb H}P^\infty )/\text{torsion})$ becomes non-abelian.
\endabstract

\endtopmatter

\document

\baselineskip =18.70pt

\heading {1. Introduction}
\endheading

Let $X$ be a connected $CW$-space and let $X^{(n)}$ denote the $n$th Postnikov approximation of $X$. The Postnikov approximations can be obtained by attaching cells to the original spaces to kill off their homotopy groups in dimensions greater than $n$. We recall that two $CW$-spaces $X$ and $Y$ are said to have the {\it same $n$-type} if  $X^{(n)}$ and $Y^{(n)}$ are homotopy equivalent.
In 1957 Adams [1] gave the example of two different homotopy types, namely $X$ and $Y$, whose Postnikov approximations are homotopy equivalent for each $n$ (see also [10]).

Let $\Sigma$ and $\Omega$ be the suspension and loop functors in the pointed homotopy category, respectively. It is
well known that the functors $\Sigma$ and $\Omega$ are examples of adjoint functors. Moreover, co-H-spaces and H-spaces are important objects of research in homotopy theory and they are the dual notions in the sense of Eckmann and Hilton (for a survey of co-H-spaces and H-spaces, see [4]).
Let $\text{Aut}(X)$ be the group of homotopy classes of self-homotopy equivalences of a space $X$, and
let $SNT(X)$ denote the set of all homotopy types $[Y]$ of spaces $Y$ such that $X$ and $Y$ have the same $n$-type for each non-negative integer $n$. We refer to Arkowitz's paper [3] and Rutter's book [21] for a survey of the vast literature about $\text{Aut}(X)$ and related topics.
In 1992 Harper and Roitberg [11] developed the connection between the set of phantom maps from $X$ to $Y$ and the set $SNT(V)$ of homotopy types of spaces $W$ having the same $n$-type as $V$ for all $n$, where $V = X \times \Omega Y$ or $Y \vee \Sigma X$.
In 1976 Wilkerson [25, theorem I] classified, up to homotopy, $CW$-spaces having the same $n$-type for all $n$ as follows: For a connected $CW$-complex $X$, there is a one-to-one correspondence between the pointed sets, $SNT(X)$ and  ${\displaystyle \lim_{\leftarrow n}}^1 \{ \text{Aut} (X^{(n)} ) \}$,
where ${\displaystyle \lim_{\leftarrow n}}^1 (-)$ is the first derived limit of groups (not necessarily abelian) in the
sense of Bousfield-Kan [7, p. 251]. Thus, if $X$ is a space of finite type, then the torsion subgroup in $\pi_* (X^{(n)})$ can be ignored in the ${\displaystyle \lim_{\leftarrow n}}^1$-calculation.

It is well known that the set of all homotopy types of spaces having the same $n$-type as the suspension of the Eilenberg-MacLane space $K( {\Bbb Z}, 2s+1))$ consists of a single homotopy type of itself; that is, $SNT(\Sigma K( {\Bbb Z}, 2s+1)) = *$. It is natural to ask for the first interesting case, i.e., the Eilenberg-MacLane space
$K( {\Bbb Z}, 2)( ={\Bbb C}P^{\infty})$ of type $( {\Bbb Z}, 2)$: Is $SNT(\Sigma {\Bbb C}P^{\infty}) = *$ or not (see [18, p. 287])?
The answer to this question was given in [13] (for a general case, see [14]). In the same lode as in the case of $\Sigma X$, where $X = {\Bbb C}P^{\infty}$ or $K( {\Bbb Z}, 2n)$, it would be interesting to know the set of all same $n$-type structure for the suspension of the infinite quaternionic projective space $\Sigma {\Bbb H}P^{\infty}$: What is $SNT(\Sigma {\Bbb H}P^{\infty})$? The major goal of this article is to give an answer to this query (see Theorem 3.1).

We note that the infinite quaternionic projective space ${\Bbb H}P^{\infty}$ is the classifying space of a compact Lie group $Sp(1)(=SU(2))$, and it has the rational homotopy type of the Eilenberg-MacLane space $K({\Bbb Z}, 4)$. Historically, McGibbon and M{\o}ller [17] showed that if $G$ is a connected compact Lie group, then its classifying space usually has an uncountable $SNT(BG)$ except for several cases, and gave an excellent set of examples. Furthermore, in [16] the classical projective
$n$-spaces (real, complex and quaternionic) were studied in terms of their self-maps from a homotopy point of view.

In this article, all spaces will be based and have the based homotopy type of a connected $CW$-complex. All maps and homotopies preserve base points. We will not distinguish notationally between a map and its homotopy class. This paper is organized as follows: In Section 2 we create iterated commutators $[\rho_{i_k},[\rho_{i_{k-1}},\ldots,[\rho_{i_{1}}, \rho_{i_2}]\ldots]]$ of self-maps $\rho_{i_j} : \Sigma \Bbb HP^\infty \rightarrow \Sigma \Bbb HP^\infty , j= 1,2,\ldots,k$ on the suspension structure of the infinite quaternionic projective space, and then show that the image of homomorphisms induced by the adjoint of iterated commutators $[\rho_{i_k},[\rho_{i_{k-1}},\ldots,[\rho_{i_{1}}, \rho_{i_2}] \ldots]]$ is both primitive and decomposable. In Section 3 we essentially show that the set of all homotopy types of spaces having the same $n$-type as the suspension of the infinite quaternionic projective space is the one element set consisting of a single homotopy type.
In Appendix, we show that the group $\text{Aut} (\pi_{\leq n} (\Sigma {\Bbb H}P^\infty )/\text{torsion})$ of automorphisms of the graded $\Bbb Z$-module, $\pi_{\leq n} (\Sigma {\Bbb H}P^\infty )$ ($= \pi_{*} (\Sigma {\Bbb H}P^\infty )^{(n)}$), preserving the Whitehead product pairings on the suspension of the infinite quaternionic projective space is finite for $n \leq 9$, and infinite for $n \geq 13$, and that $\text{Aut}(\pi_{*} (\Sigma {\Bbb H}P^\infty )/\text{torsion})$ becomes non-abelian. We obtain a similar result in the case of the infinite complex projective space. This can be extended by the use of rational homotopy theory. In particular, the techniques of localizations [15] and rationalizations of a space could be applicable.

\bigskip

\heading {2. Iterated commutators on the suspension and loop structures}
\endheading

\bigskip

For a notational convenience, we always write $Z$ for the infinite quaternionic projective space ${\Bbb H}P^{\infty}$. Let $\hat \rho_1 : Z \rightarrow \Omega \Sigma Z$ be the canonical inclusion map, and inductively define (cf. [20, p. 153])
$$
\hat \rho_{n+1} : Z @>{\bar \Delta}>> Z \wedge Z @>{\hat \rho_{1} \wedge \hat
\rho_{n}}>> \Omega \Sigma Z \wedge \Omega \Sigma Z @>{C}>> \Omega \Sigma Z ,
$$
where $\bar \Delta$ is the reduced diagonal map (i.e., the composite of the diagonal $\Delta : Z \rightarrow Z \times Z$ with the projection $\pi : Z \times Z \rightarrow Z \wedge Z$ onto the smash product) and $C : \Omega \Sigma Z \wedge \Omega \Sigma Z \rightarrow \Omega \Sigma Z$ is the commutator map with respect to the loop operation; that is,
$$
\hat \rho_{n+1}(z) = C(\hat \rho_{1} (z), \hat \rho_{n} (z) ) = \hat \rho_{1}(z) \cdot \hat \rho_{n}(z) \cdot (\hat \rho_{1}(z))^{-1} \cdot (\hat \rho_{n}(z))^{-1}.
$$
Here the multiplication is the loop multiplication and the inverse means the loop inverse $l : \Omega \Sigma Z \rightarrow \Omega \Sigma Z$. The above definition does make sense because there are infinitely many non-zero cohomology cup products in $Z$, so that it has the infinite Lusternik-Schnirelmann category [24, Chapter X]. Moreover, there is a result of Arkowitz and Curjel [5, Theorem 5] which asserts that the $n$-fold commutator is of finite order if and only if all $n$-fold cup products of any positive dimensional rational cohomology classes of a space vanish (see also [2]).

By using the Serre spectral sequence of a fibration, we can see that $Z$ has a $CW$-decomposition as follows:
$$
Z = S^{4}  \cup_{\gamma_1} e^{8}   \cup_{\gamma_2} e^{12}  \cup_{\gamma_3} \cdots
 \cup_{\gamma_{n-1}}e^{4n} \cup_{\gamma_n}e^{4(n+1)} \cup_{\gamma_{n+1}}  \cdots ,
$$
where the $\gamma_n$ are the attaching maps as Hopf fibrations $\gamma_n : S^{4n+3} \rightarrow {\Bbb H}P^{n}$ with fiber $S^3$ for each positive integer $n$. Even though the formations of $CW$-decompositions of the infinite complex and quaternionic projective spaces are very much alike, there is an important difference between them in that the infinite complex projective space ${\Bbb C}P^{\infty}$ is an H-space, while $Z (= {\Bbb H}P^{\infty})$ is not an H-space but a rational H-space; that is,
$$
Z ~ \simeq_{\Bbb Q}~ K(\Bbb Z,4)~ \simeq_{\Bbb Q} ~ S^4 \vee S^8 \vee \cdots \vee S^{4n} \vee \cdots ,
$$
where $\simeq_{\Bbb Q}$ is a rational homotopy equivalence. Thus, by using the Hilton-Milnor theorem [24, p. 511], we can find a rationally non-trivial indecomposable generators, say, $\hat x_n : S^{4n} \rightarrow \Omega \Sigma Z$, in the homotopy groups $\pi_{4n}(\Omega \Sigma Z)$ for each $n = 1,2,3,\ldots$. Indeed, it is possible for us to do so because there is a one-to-one correspondence between the integral generators modulo torsions and rational generators.

Now we take the self-maps $\rho_n : \Sigma Z \rightarrow \Sigma Z$ and maps $x_n : S^{4n+1} \rightarrow \Sigma Z$ by using the adjoint relations of $\hat \rho_{n} : Z \rightarrow \Omega \Sigma Z$ and $\hat x_n : S^{4n} \rightarrow \Omega \Sigma Z$, respectively, for $n = 1,2,3, \ldots$.

\bigskip

\remark{Remark 2.1} Let $[\rho_{i_1}, \rho_{i_2}] : \Sigma Z \rightarrow \Sigma Z$ be the commutator of self-maps, $\rho_{i_1}$ and $\rho_{i_2}$, on the suspension structure; that is,
$[\rho_{i_1}, \rho_{i_2}] = \rho_{i_1} + \rho_{i_2} - \rho_{i_1} - \rho_{i_2}$,
where the operations are the suspension additions on $\Sigma Z$.
We construct self-maps of $\Sigma Z$ by $I + [\rho_{i_k},[\rho_{i_{k-1}},\ldots,[\rho_{i_{1}}, \rho_{i_2}] \ldots]]$, where $I$ is the identity map on $\Sigma Z$ and $[\rho_{i_k},[\rho_{i_{k-1}},\ldots,[\rho_{i_{1}}, \rho_{i_2}]\ldots]]$ is the iterated commutator of self-maps $\rho_{i_j} : \Sigma Z \rightarrow \Sigma Z$ on the suspension structure for $j= 1,2,\ldots,k$. These self-maps play a key role for answering to the given same $n$-type question in that those self-maps
$I + [\rho_{i_k},[\rho_{i_{k-1}},\ldots,[\rho_{i_{1}}, \rho_{i_2}] \ldots]]$ are self-homotopy equivalences due to the Whitehead theorem. We can see that the above self-homotopy equivalences are related to the Lie algebra homomorphisms $\theta : \pi_* (\Sigma K(\Bbb Q,4)) \rightarrow \pi_* (\Sigma K(\Bbb Q,4))$ such that $\theta (z_i ) = z_i + \Sigma_{s \geq 2} Q_s (z_i )$, where $Q_s (z_i )$ is the brackets of lower dimensional elements of length $s$.
\endremark

\bigskip

We recall that
$$
\tilde H_* ( Z ; {\Bbb Q}) \cong {\Bbb Q} \{ b_1 , b_2 , \ldots, b_n , \ldots \},
$$
as a graded $\Bbb Q$-module, where $b_n $ is the standard generator of $H_{4n} (Z ; {\Bbb Q})$ for each $n = 1,2,3, \ldots$.  The Bott-Samelson theorem [6] says that the Pontryagin algebra $H_* (\Omega \Sigma Z ;\Bbb Q )$ is isomorphic to the tensor algebra $TH_* (Z ; \Bbb Q )$; that is, the rational homology of $\Omega \Sigma Z $ is the tensor algebra $T [b_1, b_2,  \ldots, b_n , \ldots ]$ generated by  $\{ b_1 , b_2 , \ldots, b_n , \ldots \}$.
From the Serre spectral sequence of a fibration, we have an isomorphism $H^* (Z;\Bbb Q) \cong \Bbb Q[\lambda]$ as algebras. Here $\Bbb Q[\lambda]$ is the polynomial algebra over $\Bbb Q$ generated by $\lambda$ of degree $4$; that is, $\lambda$ is a generator of $H^{4}(Z;\Bbb Q)$ so that the Kronecker index $<\lambda^i ,b_j>=\delta_{ij}$.

Let $i_1 : X \rightarrow X \times X$ and $i_2 : X \rightarrow X \times X$ be the inclusions defined by $i_1 (x) = (x, x_0 )$ and $i_2 (x) = (x_0 , x)$, respectively, where $x_0$ is the base point of $X$. Recall that an element $z \in H_* (X)$ is said to be {\it primitive} if and only if $\Delta_* (z) = {i_1}_* (z) + {i_2}_* (z) = z \otimes 1 + 1 \otimes z$ in homology, where $\Delta : X \rightarrow X \times X$ is the diagonal map.

Since $\Sigma Z$ is a co-H-space (actually, cogroup-like), we have $(\rho_m + \rho_n )_* = {\rho_m}_* + {\rho_n}_*$ and $(\rho_m + \rho_n )^* = {\rho_m}^* + {\rho_n}^*$ in homology and cohomology levels, respectively, where the additions in the parentheses are suspension sums. However, in general those results do not hold in the case of multiplications $\hat \rho_m \cdot \hat \rho_n$ on H-spaces, where `$\cdot$' is the loop multiplication. Fortunately, the above equalities could be guaranteed for the action of loop multiplications on the {\it primitive} elements. This raises a basic question: What types of elements could be primitive as well as decomposable? We give an answer to this question as follows:

\bigskip

\proclaim{Theorem 2.2} Let $[\hat \rho_{i_k},[\hat\rho_{i_{k-1}},\ldots,[\hat\rho_{i_{1}}, \hat\rho_{i_2}] \ldots]]$ be an iterated commutator in the loop structure. Then $[\hat \rho_{i_k},[\hat\rho_{i_{k-1}},\ldots,[\hat\rho_{i_{1}}, \hat\rho_{i_2}] \ldots]]_*(b_{i_1 + i_2 +\cdots + i_k})$ is both primitive and decomposable for each $i_j, j=1,2,\ldots,k$, where $b_{i_1 + i_2 +\cdots + i_k}$ is a generator in the tensor algebra $T [b_1, b_2,  \ldots, b_n , \ldots ]$.
\endproclaim

\bigskip

\demo{Proof} Let $X_p$ denote the $p$-skeleton of a $CW$-complex $X$. In order to prove this theorem, we need
the following lemmas which will be also used in Section 3:

\bigskip

\proclaim{Lemma 2.3} Let $X$ be a finite type $CW$-complex with base point $x_0$ as the zero skeleton and let $f,g : X \rightarrow \Omega Y$ be base point preserving maps with $f \vert_{X_p} \simeq *$ and $g \vert_{X_q} \simeq *$. Then the restriction of the commutator $[f, g] : X \rightarrow \Omega Y$ to the $(p+q)$-skeleton of $X$ is inessential.
\endproclaim

\demo{Proof} We can factor the commutator as follows:
$$
\CD
X @>[f,g]>> \Omega Y \\
@V{\Delta}VV  @V=VV \\
X \times X @>{C(f,g)}>>  \Omega Y ,
\endCD
$$
where $C(f,g)$ sends $(x,y) \mapsto f(x)\cdot g(y)\cdot (f(x))^{-1}\cdot (g(y))^{-1}$. Here the multiplication means the loop multiplication. We can give $X\times X$ the product $CW$-structure with
$$
(X \times X)_n = \bigcup_{j=0}^{n} X_j \times X_{n-j}.
$$
We note that if $j \leq p$ or $n-j \leq q$, then the restrictions $C(f,g)\mid_{X_j \times X_{n-j}}$ are null-homotopic, say,
null-homotopies $F_{j, n-j} : C(f,g)\mid_{X_j \times X_{n-j}} \simeq *$. Using the homotopy extension property, we are able to
extend those null-homotopies to the $(p+q)$-skeleton $(X \times X)_{p+q}$. Thus the cellular approximation theorem shows the proof. \hfill\qquad$\square$
\enddemo

\bigskip

\proclaim{Lemma 2.4} Let $f:X \rightarrow Y$ be a map such that $f\vert_{X_p}\, = \,*$. Then the image of $f_*$ in $H_{p+1} (Y)$ is primitive.
\endproclaim

\demo{Proof} Since $f\vert_{X_p} \simeq *$, the map $f : X \rightarrow Y$ can be factored as follows:
$$
\CD
X @>f>> Y \\
@VqVV @V=VV \\
X /X_p @>{\bar f}>> Y
\endCD
$$
However, since $X/X_p$ is $p$-connected, by the Hurewicz isomorphism theorem, every class in $H_{p+1}
(X/X_p )$ is spherical, and thus primitive. We therefore know that the image of $\bar f_*$ (and thus the one of $f_*$) lies in the set of primitives $PH_{p+1}(Y)$ in $H_{p+1}(Y)$.  \hfill\qquad$\square$
\enddemo

\bigskip

\remark{Remark 2.5} Let $\widehat {[\rho_{i_1}, \rho_{i_2}]} : Z \rightarrow \Omega \Sigma Z$ be the
adjoint of the commutator $[\rho_{i_1}, \rho_{i_2}] : \Sigma Z \rightarrow \Sigma Z$. Since the map \,\, $\kappa = ~~\widehat{}~~ : [\Sigma Z, \Sigma Z] \rightarrow [Z, \Omega \Sigma Z]$, defined by
$$
(\kappa{\varphi})(z)(t) = \hat{\varphi}(z)(t) = \varphi (\prec z,t \succ ),
$$
for $\varphi \in [\Sigma Z, \Sigma Z], z \in Z, t \in I$ and $\prec z,t \succ \in \Sigma Z$, is a group isomorphism, it follows that $\widehat {[\rho_{i_1}, \rho_{i_2}]} = [\hat \rho_{i_1}, \hat \rho_{i_2}] : Z \longrightarrow \Omega \Sigma Z$, the commutator of $\hat \rho_{i_1}$ and $\hat \rho_{i_2}$ in $[Z, \Omega \Sigma Z]$. Similarly, we have
$$
\kappa({[\rho_{i_k},[\rho_{i_{k-1}},\ldots,[\rho_{i_{1}}, \rho_{i_2}] \ldots]]}) =
[\hat \rho_{i_k},[\hat\rho_{i_{k-1}},\ldots,[\hat\rho_{i_{1}}, \hat\rho_{i_2}] \ldots]].
$$
\endremark

\bigskip

Let $Z_j$ denote the $j$-skeleton of $Z$. Then the following diagram
$$
\CD
Z_{4m+4n} @>i>> Z  \\
@VVpV   @VV{[\hat \rho_m , \hat \rho_n ]}V \\
S^{4m+4n} @>{<\hat x_m , \hat x_n > }>> \Omega \Sigma Z
\endCD \tag 2-1
$$
is commutative up to homotopy.
Here $i$ is the inclusion, $p$ is the projection and $<\hat x_m , \hat x_n >$ is the Samelson product.
Moreover, the map $\hat \rho_{n} : Z \rightarrow \Omega \Sigma Z$ can be factored as
$$
Z @>q>> Z/Z_{4n-1} @>g_n>> \Omega \Sigma Z
$$
such that the restriction to the bottom sphere of the map
$$
g_n : Z/Z_{4n-1}(= S^{4n} \cup_{\gamma_n} e^{4n+4} \cup_{\gamma_{n+1}} \cdots) ~ \longrightarrow ~\Omega \Sigma Z
$$
coincides with a map $S^{4n} \rightarrow \Omega \Sigma Z$; that is $\hat \rho_{n} |_{Z_{4n-1}} \simeq *$, where $q : Z \rightarrow Z/Z_{4n-1}$ is the projection. Indeed, the proof of this statement is analogous to the proof of the Morisugi's result [20, Theorem 1.4] in the case of the infinite complex projective space.

\bigskip

\proclaim{Lemma 2.6} $[\rho_{i_k},[\rho_{i_{k-1}},\ldots,[\rho_{i_{1}}, \rho_{i_2}] \ldots]]|_{(\Sigma Z)_{4(i_1 + \cdots +i_k )}}$ is inessential for each $\rho_{i_j} : \Sigma Z \rightarrow \Sigma Z, j =1, 2, \ldots,k$.
\endproclaim

\demo{Proof}
Since $\widehat{[\rho_{i_1} , \rho_{i_2} ]} = [ \hat{\rho_{i_1}}, \hat{\rho_{i_2}}]$ and $\hat \rho_{i_j} |_{Z_{4i_j -1}} \simeq *$ for $j = 1,2$, by Lemma 2.3, we see that the commutator $[ \hat{\rho_{i_1}},
\hat{\rho_{i_2}}]$ restricts to the trivial map on $Z_{4(i_1 + i_2) -2}$. Since $Z$ has no cells in
dimensions $4(i_1 + i_2) -3, 4(i_1 + i_2) -2$ and $4(i_1 + i_2) -1$, we have
$$
Z_{4(i_1 + i_2) -3} = Z_{4(i_1 + i_2) -2} = Z_{4(i_1 + i_2) -1}.
$$
By taking the adjointness, we reach the proof in the case of two-fold commutators.

For the general case, we suppose that $[\hat\rho_{i_{k-1}},\ldots,[\hat\rho_{i_{1}}, \hat\rho_{i_2}] \ldots]|_{Z_{4(i_1 + \cdots +i_{k-1} )-1}}$ is inessential. Then since $\hat \rho_{i_k} |_{Z_{4i_k -1}} \simeq *$, the similar argument as described above shows that $[\hat \rho_{i_k},[\hat\rho_{i_{k-1}},\ldots,[\hat\rho_{i_{1}}, \hat\rho_{i_2}] \ldots]]|_{Z_{4(i_1 + \cdots +i_k )-1}}$ is nullhomotopic by induction on $k$, and the adjointness shows the proof of this lemma. \hfill\qquad$\square$
\enddemo

\bigskip

What is the relationship between the commutators of the loop and suspension structures? The following lemma provides the answer to this query.

\bigskip

\proclaim{Lemma 2.7} $(\Omega [\rho_{i_k},[\rho_{i_{k-1}},\ldots,[\rho_{i_{1}}, \rho_{i_2}] \ldots]]) \circ E
= [\hat \rho_{i_k},[\hat\rho_{i_{k-1}},\ldots,[\hat\rho_{i_{1}}, \hat\rho_{i_2}] \ldots]]$, where $E : Z \rightarrow
\Omega \Sigma Z$ is the canonical inclusion map.
\endproclaim

\bigskip

\demo{Proof} We can see that the embedding $E : Z \rightarrow \Omega \Sigma Z$ sends $z \mapsto E(z) : I \rightarrow \Sigma Z$, and  $E(z)(t) = \prec z,t \succ ~~\in~~ \Sigma Z$. From the composition
$$
(\Omega [\rho_{i_k},[\rho_{i_{k-1}},\ldots,[\rho_{i_{1}}, \rho_{i_2}] \ldots]]) \circ E : Z \longrightarrow \Omega \Sigma Z,
$$
of maps, the map $((\Omega [\rho_{i_k},[\rho_{i_{k-1}},\ldots,[\rho_{i_{1}}, \rho_{i_2}] \ldots]]) \circ E)(z) : I \rightarrow \Sigma Z$ sends $t \in I$ to
$$
\split
((\Omega [\rho_{i_k},[\rho_{i_{k-1}},\ldots,[\rho_{i_{1}}, \rho_{i_2}] \ldots]]) \circ E)(z)(t) &= [\rho_{i_k},[\rho_{i_{k-1}},\ldots,[\rho_{i_{1}}, \rho_{i_2}] \ldots]] (E(z)(t)) \\
&= [\rho_{i_k},[\rho_{i_{k-1}},\ldots,[\rho_{i_{1}}, \rho_{i_2}] \ldots]](\prec z,t \succ),
\endsplit
$$
where $z \in Z$.

Let $\kappa : [\Sigma Z, \Sigma Z] \rightarrow [Z, \Omega \Sigma Z]$ be the adjoint isomorphism as in Remark 2.5. Then
$\kappa([\rho_{i_k},[\rho_{i_{k-1}},\ldots,[\rho_{i_{1}}, \rho_{i_2}] \ldots]]) : Z \rightarrow \Omega \Sigma Z$ sends $z \in Z$ to
$$
\kappa({[\rho_{i_k},[\rho_{i_{k-1}},\ldots,[\rho_{i_{1}}, \rho_{i_2}] \ldots]]})(z) : I \longrightarrow \Sigma Z.
$$
Here, $(\kappa({[\rho_{i_k},[\rho_{i_{k-1}},\ldots,[\rho_{i_{1}}, \rho_{i_2}] \ldots]]})(z))(t) = [\rho_{i_k},[\rho_{i_{k-1}},\ldots,[\rho_{i_{1}}, \rho_{i_2}] \ldots]](\prec z,t \succ)$. This completes the proof. \hfill\qquad$\square$
\enddemo

\bigskip

The above lemmas show that $[\hat \rho_{i_k},[\hat\rho_{i_{k-1}},\ldots,[\hat\rho_{i_{1}}, \hat\rho_{i_2}] \ldots]]_* (b_{i_1 + \cdots +i_k})$ is primitive. In order to prove that it is decomposable, we use the homology suspension $\sigma : H_{4n} (\Omega \Sigma Z  ) \longrightarrow H_{4n+1} (\Sigma Z )$. We note that if $z \in H_{4n+1} (\Sigma Z )$ is spherical, then $z$ is in the image of $\sigma$.
This follows since $\sigma$ is natural and thus the following diagram commutes:
$$
\CD
H_{4n} (\Omega S^{4n+1}) @>{\sigma}>\cong> H_{4n+1}(S^{4n+1}) \\
@V{\Omega f_*}VV @V{f_*}VV \\
H_{4n} (\Omega \Sigma Z  ) @>{\sigma}>>  H_{4n+1}(\Sigma Z ).
\endCD
$$
In particular, an element $\lambda \beta_n$ of the image $\sigma (H_{4n}(\Omega \Sigma Z  ))$ in $H_{4n+1}(\Sigma Z )$ is spherical  for each $n$, and $\sigma$ takes all products (decomposable terms)
to zero (see [24, p. 383]), where $\lambda \neq 0$ and $\beta_n$ is a generator of $H_{4n+1}(\Sigma K )$. Therefore, we have
$\sigma (b_{n}) \,= \, \lambda \beta_{n}$ for some $\lambda \neq 0$,
since $b_{n}$ is indecomposable. On the other hand, the homomorphism $[\rho_{i_k},[\rho_{i_{k-1}},\ldots, [\rho_{i_{1}}, \rho_{i_2}] \ldots]]_*$ induced by the iterated commutator is trivial  in homology (except in degree $0$).
We now observe that the embedding $E (= \hat \rho_1 ) : Z \rightarrow \Omega \Sigma Z$ has a left homotopy inverse rationally, and thus, by using the universal coefficient theorem for homology, we have a monomorphism $E_* : H_* (Z; \Bbb Q) \rightarrow H_* (\Omega \Sigma Z; \Bbb Q)$ so that $E_* (b_n ) = b_n$ for each $n$.
Using the following commutative diagram
$$
\CD
H_{4(i_1 + \cdots + i_k )} (\Omega \Sigma Z ) @>\sigma>> H_{4(i_1 + \cdots + i_k )+1}(\Sigma Z ) \\
@V{\Omega [\rho_{i_k},[\rho_{i_{k-1}},\ldots, [\rho_{i_{1}}, \rho_{i_2}] \ldots]]_{*}}VV
@VV{[\rho_{i_k},[\rho_{i_{k-1}},\ldots, [\rho_{i_{1}}, \rho_{i_2}] \ldots]]_*}V \\
H_{4(i_1 + \cdots + i_k )} (\Omega \Sigma Z  ) @>{\sigma}>>  H_{4(i_1 + \cdots + i_k )+1} (\Sigma Z)
\endCD
$$
together with Lemma 2.7, we can see that  $[\hat \rho_{i_k},[\hat\rho_{i_{k-1}},\ldots,[\hat\rho_{i_{1}}, \hat\rho_{i_2}] \ldots]]_* (b_{i_1 + \cdots +i_k})$ is decomposable as required. \hfill\qquad$\square$
\enddemo

\bigskip

\heading {3. Triviality of $SNT(\Sigma Z)$}
\endheading

\bigskip

The main purpose of this section is to prove the following theorem:

\proclaim{Theorem 3.1} $SNT(\Sigma Z) = *$.
\endproclaim

\demo{Proof} By using the results of previous section, we prove this theorem as follows:

Let $\text{Aut}( \pi_{\leq n} (X))$ denote the group of automorphisms of the non-negatively graded group, $\pi_{\leq n} (X)$, preserving the Whitehead product pairings. In 1992 McGibbon and M{\o}ller [18, Theorem 1] proved the following theorem:

\bigskip

\proclaim{Theorem 3.2} Let $X$ be a $1$-connected space with finite
type over some subring of the rationals. Assume that $X$ has the rational homotopy type of a bouquet of spheres. Then the following three conditions
are equivalent: \roster
\item $SNT(X) = *$;
\item the map
$\text{Aut} (X) @>{f \mapsto f^{(n)}}>> \text{Aut}(X^{(n)})$ has a finite cokernel for each $n$; and
\item the map
$\text{Aut} (X) @>{f \mapsto f_\sharp}>> \text{Aut}( \pi_{\leq n} (X))$ has a finite cokernel for each $n$.
\endroster
\endproclaim

\bigskip

Recall that the Samelson product
$$
< ~,~ > \, \, : \pi_{m-1}(\Omega X) \times \pi_{n-1}(\Omega X)
\rightarrow \pi_{m+n-2} (\Omega X)
$$
is defined by the commutative diagram
$$
\CD
\pi_{m}(X) \times \pi_{n}(X) @>{[~,~]}>> \pi_{m+n-1} (X) \\
@V{\tau \times \tau}V{\cong}V  @V{(-1)^{m-1} \tau}V{\cong}V \\
\pi_{m-1}(\Omega X) \times \pi_{n-1}(\Omega X) @>{<~,~>}>>
\pi_{m+n-2} (\Omega X),
\endCD
$$
where $[~,~]$ means the Whitehead product and the vertical arrows are given by the adjointness relation $\pi_m (X) = [\Sigma S^{m-1}
,X] \cong [S^{m-1} , \Omega X] = \pi_{m-1}(\Omega X)$.

We observe that the total rational homotopy group $\hat \Cal L = \pi_{*} (\Omega \Sigma Z ) \otimes \Bbb Q$ of $\Omega \Sigma Z$ becomes a graded Lie algebra over $\Bbb Q$ with the Lie bracket $< ~,~>$ given by the Samelson product which is called the {\it rational homotopy Lie algebra} [23] of $\Sigma Z$. Let $\hat \Cal L_{\leq n}$ denote the subalgebra of $\hat \Cal L$ generated by all free algebra generators of degree less than or equal to $4n$. That is, $\hat \Cal L_{\leq n} =  \pi_{\leq 4n} (\Omega \Sigma Z ) \otimes \Bbb Q$ with generators $\hat \chi_1 , \hat \chi_2 , \ldots ,\hat \chi_n $, where $\hat \chi_i : S^{4i} @>>> \Omega \Sigma Z_{\Bbb Q}$ is the composition $r \circ \hat x_i$ of the rationally non-trivial indecomposable element $\hat x_i : S^{4i} @>>> \Omega \Sigma Z$ of $\pi_{4i} (\Omega \Sigma Z )$, $i = 1,2, \ldots, n$ with the rationalizing map $r : \Omega \Sigma Z @>>> \Omega \Sigma Z_{\Bbb Q}$. Similarly, we can also consider $\Cal L_{\leq n} =  \pi_{\leq 4n+1} (\Sigma Z ) \otimes \Bbb Q$ with generators $\chi_1 , \chi_2 , \ldots , \chi_n$. Here $\chi_i : S^{4i+1} @>>> \Sigma Z_{\Bbb Q}$ is the composite $r \circ x_i$ of the rationally non-trivial indecomposable element $x_i : S^{4i+1} @>>> \Sigma Z$ of $\pi_{4i+1} (\Sigma Z )$ with the rationalizing map $r : \Sigma Z @>>> \Sigma Z_{\Bbb Q}$.
Furthermore, $\Cal L_{\leq n}$ given by the Whitehead product $[~,~]$ has the graded quasi-Lie algebra structure which is called the {\it Whitehead algebra}.

\bigskip

\remark{Remark 3.3} Let $h : \pi_* (\Omega \Sigma Z) \rightarrow H_* (\Omega \Sigma Z ; \Bbb Q)$ be the Hurewicz homomorphism. Then by using the commutative diagram (2-1) in Section 2, we can obtain the similar results from Proposition 3.1 and Lemma 3.2 in [13] for the infinite quaternionic projective space as follows:
\roster
\item $[\hat \rho_m , \hat \rho_n ]_*
(b_{m+n}) = [h(\hat x_m ) , h(\hat x_n )] ( = h <\hat x_m , \hat x_n >)$ in rational homology (see [8, p. 141] for the second equality). Here $b_{m+n}$ is the standard generator of $H_{4m+4n} (Z ; \Bbb Q)$, $\hat x_m$ and $\hat x_n$ are rationally non-trivial indecomposable elements in homotopy groups, and $[h(\hat x_m ) , h(\hat x_n )] = h(\hat x_m )h(\hat x_n ) - (-1)^{\vert h(\hat x_m )\vert \vert h(\hat x_n )\vert} h(\hat x_n )h(\hat x_m )$.
\item Let $x_{n}$ be the rationally non-trivial indecomposable element in $\pi_{4n+1} (\Sigma Z)$ for $n = i_1 + i_2 + \cdots + i_k$. Then
$(\text{I} + [\rho_{i_k},[\rho_{i_{k-1}},\ldots,[\rho_{i_{1}}, \rho_{i_2}]\ldots]])_\sharp (x_n ) \,=\, x_n + [\rho_{i_k},[\rho_{i_{k-1}},\ldots,[\rho_{i_{1}}, \rho_{i_2}]\ldots]]_\sharp (x_n )$,
where the first $+$ is the suspension addition on $\Sigma Z$, while the second $+$ refers to the addition of homotopy groups.
\endroster
\endremark

\bigskip

\remark{Remark 3.4} The famous Cartan-Serre theorem [9, p. 231] asserts that the Hurewicz homomorphism $h : \pi_* (\Omega \Sigma Z) \rightarrow H_* (\Omega \Sigma Z ; \Bbb Q)$ induces an isomorphism
$$
\pi_* (\Omega \Sigma Z)\otimes \Bbb Q \cong PH_* (\Omega \Sigma Z; \Bbb Q),
$$
where $PH_* (\Omega \Sigma Z; \Bbb Q)$ means the primitive subspace of $H_* (\Omega \Sigma Z; \Bbb Q)$ (see also [19]).
Moreover, we have $h(\hat x_n ) = \alpha b_n + \text{decomposables}, \alpha \neq 0$ for each $n = 1,2,3,\ldots$. Here $\hat x_n$ is the rationally non-trivial indecomposable element of $\pi_{4n}(\Omega \Sigma Z)$, $b_n ( = E_* (b_n ))$ is the generator of $H_{4n}(\Omega \Sigma Z ; \Bbb Q)$, and the decomposable parts consist of products of generators with dimensions less than $4n$ in the tensor algebra.
\endremark

\bigskip

\proclaim{Lemma 3.5} The iterated Whitehead products $[x_{i_k},[x_{i_{k-1}},\ldots,[x_{i_{1}}, x_{i_2}]\ldots]]$ in homotopy groups $\pi_{4(i_1 + \cdots +i_k)+1}(\Sigma Z )$ are rationally non-trivial, where the $x_{i_j}$ are indecomposable generators in dimension
$4i_j +1 , j = 1,\ldots,k$.
\endproclaim

\demo{Proof} We first prove the result in the case of two-fold Whitehead product. Suppose that
$[x_{i_1}, x_{i_2}]$ has a finite order in $\pi_{4(i_1 +i_2)+1}(S^{4i_1 +1} \vee S^{4i_2 +1})$, then,
from the cofibration
$$
S^{4(i_1 +i_2)+1} @>>> S^{4i_1+1} \vee S^{4i_2 +1} @>>>  S^{4i_1+1} \times S^{4i_2+1} @>>> S^{4(i_1 +i_2)+2} @>>> \cdots,
$$
we have
$$
S^{4i_1 +1} \times S^{4i_2 +1} \simeq (S^{4i_1 +1} \vee S^{4i_2 +1}) \cup_{[x_{i_1}, x_{i_2}]} e^{4(i_1 +i_2)+2},
$$
where $[x_{i_1}, x_{i_2}] :  S^{4(i_1 +i_2)+1} \rightarrow S^{4i_1+1} \vee S^{4i_2 +1}$ is the attaching map.
Since the rationalized Whitehead product $[x_{i_1}, x_{i_2}]_{\Bbb Q}$ is trivial by assumption, we obtain
$$
S^{4i_1 +1} \times S^{4i_2 +1} \simeq_{\Bbb Q} S^{4i_1 +1} \vee S^{4i_2 +1} \vee S^{4(i_1 +i_2)+2} \tag 3-1
$$
which is a contradiction by using the cohomology cup product argument, where `$\simeq_{\Bbb Q}$' is a rational homotopy equivalence. Indeed, there exist non-zero cup products in $H^* (S^{4i_1 + 1} \times S^{4i_2 + 1})$, while all cup products on the right hand side of (3-1) are always trivial.

For the general case, we suppose that the result holds for the $(k-1)$-fold Whitehead products, i.e., $[x_{i_{k-1}},\ldots,[x_{i_{1}}, x_{i_2}]\ldots]$ is a rationally non-trivial in $\pi_{4(i_1 + \cdots +i_{k-1})+1}(\Sigma Z )$. By the Hilton's formula [12], we can consider the basic product  $[x_{i_{k-1}},\ldots,[x_{i_{1}}, x_{i_2}]\ldots]$ as a rational generator of
$$
\pi_{4(i_1 + \cdots +i_{k-1})+1}(S^{4(i_1 + \cdots +i_{k-1})+1} ) \otimes \Bbb Q
\subset \pi_{4(i_1 + \cdots +i_{k-1})+1}(\Sigma Z ) \otimes \Bbb Q.
$$
We note that all rational generators of $\pi_* (\Sigma Z) \otimes \Bbb Q$ consists of basic products because the rationalization of $\Sigma Z$ has the same homotopy type of a bouquet of infinitely many odd dimensional spheres. The similar method as mentioned above shows that
$$
S^{4(i_1 + \cdots +i_{k-1})+1} \times S^{4i_k +1} \simeq (S^{4(i_1 + \cdots +i_{k-1})+1} \vee S^{4i_k +1})
\cup_{\gamma} e^{4(i_1 + \cdots +i_{k-1} + i_k)+2},
$$
where $\gamma = [x_{i_k},[x_{i_{k-1}},\ldots,[x_{i_{1}}, x_{i_2}]\ldots]] : S^{4(i_1 + \cdots +i_{k})+1} \rightarrow S^{4(i_1 + \cdots +i_{k-1})+1} \vee S^{4i_k +1}$ is the attaching map.
If this map $\gamma = [x_{i_k},[x_{i_{k-1}},\ldots,[x_{i_{1}}, x_{i_2}]\ldots]]$ has a finite order in homotopy groups, then we have a contradiction again by the same assertion of the above statement as required.
\hfill\qquad$\square$
\enddemo

\bigskip

In general, we observe that the basic products $[x_{i_k},[x_{i_{k-1}},\ldots,[x_{i_{1}}, x_{i_2}]\ldots]]$ are not integral homotopy generators but rational homotopy generators when they are rationalized.

\bigskip

\proclaim{Lemma 3.6}  For each iterated Whitehead product
$[x_{i_k},[x_{i_{k-1}},\ldots,[x_{i_{1}}, x_{i_2}]\ldots]]$ in the non-negatively graded group $\pi_* (\Sigma Z)$, there exists an iterated commutator
$[\rho_{i_k},[\rho_{i_{k-1}},\ldots,$ $[\rho_{i_{1}}, \rho_{i_2}]\ldots]]$ in the group $[\Sigma Z, \Sigma Z ]$
such that
$$
(I + [\rho_{i_k},[\rho_{i_{k-1}},\ldots,[\rho_{i_{1}}, \rho_{i_2}]\ldots]])_\sharp (x_n )
= x_n +\alpha [x_{i_k},[x_{i_{k-1}},\ldots,[x_{i_{1}}, x_{i_2}]\ldots]] ,
$$
where $\alpha \neq 0$, and $x_n$ and $x_{i_j}$ are rationally non-trivial indecomposable elements and $n = i_1 + i_2 + \cdots + i_k$.
\endproclaim

\demo{Proof} Step 1. We prove this lemma by induction on $k$. First, in order to show that
$[\rho_{i_1}, \rho_{i_2}]_\sharp (x_n ) \,= \,\alpha [ x_{i_1} , x_{i_2}]$, where $\alpha \neq 0$ and $n= i_1 + i_2$, we consider the following commutative diagram:
$$
\CD
\pi_* (\Omega \Sigma Z ) @>{\Omega[\rho_{i_1} , \rho_{i_2}]_\sharp }>> \pi_* (\Omega \Sigma Z )  \\
@VVhV  @VVhV \\
H_* (\Omega \Sigma Z ; \Bbb Q ) @>{\Omega[\rho_{i_1} , \rho_{i_2}]_* }>> H_* (\Omega \Sigma Z ; \Bbb Q ).
\endCD
$$
For the rationally non-trivial indecomposable element $\hat x_n$ in $\pi_{4n} (\Omega \Sigma Z )/\text{torsion}$,
we have
$$
\split
h \Omega[&\rho_{i_1}, \rho_{i_2}]_\sharp (\hat x_n ) = \Omega [\rho_{i_1}, \rho_{i_2}]_* h (\hat x_n ) \hskip 1.6cm (\text{commutativity}) \\
&= \Omega [\rho_{i_1}, \rho_{i_2}]_* ( \alpha b_n + \text{decomposables}) \hskip 0.8cm (\alpha \neq 0 \,\, \text{by Remark 3.4} )\\
&= \widehat {[\rho_{i_1}, \rho_{i_2}]}_* (\alpha b_n + \text{decomposables}) \hskip 1cm (E_* (b_n ) = b_n ~~ \text{and by Lemma 2.7} )\\
&= [\hat \rho_{i_1}, \hat \rho_{i_2}]_* (\alpha b_n ) + 0 = \alpha [h(\hat x_{i_1} ), h(\hat x_{i_2} )] \hskip 0.2cm (\text{by Lemma 2.6, and Remarks 2.5 and 3.3} )\\
&= \alpha h( <\hat x_{i_1} , \hat x_{i_2} > ) = h( \alpha <\hat x_{i_1} , \hat x_{i_2} > ).
\endsplit \tag 3-2
$$
It can be noticed that the above zero term is derived from the fact that
the restriction $[\hat \rho_{i_1}, \hat \rho_{i_2}]\vert_{Z_{4(i_1 +i_2 ) -1}}$ to the skeleton is inessential due to the proof of Lemma 2.6; that is, $[\hat \rho_{i_1}, \hat \rho_{i_2}]_* (b_j ) =0$ for $\text{dim} (b_j ) < 4(i_1 +i_2 )$ in rational homology.

We see that, by Lemma 3.5 and adjointness, $\hat x_n$ and $< \hat x_{i_1} , \hat x_{i_2}>$ are rationally non-trivial indecomposable and decomposable elements, respectively, in $\pi_{4n} (\Omega \Sigma Z)/\text{torsion}$, and that $h(\hat x_n )$ is spherical and thus primitive. By considering the above equation (3-2) and the Cartan-Serre theorem described in
Remark 3.4, we can conclude that $\Omega[\rho_{i_1}, \rho_{i_2}]_\sharp (\hat x_n ) \,\, = \,\,\alpha <\hat x_{i_1}, \hat x_{ i_2} >$. We also remark that the map $\pi_* (\Omega \Sigma Z)/\text{torsion} \rightarrow H_* (\Omega \Sigma Z ; \Bbb Q)$ is a monomorphism, and that $\Omega [\rho_{i_1}, \rho_{i_2}]$ is a loop map, thus it is an H-map. Moreover, there is a bijection between $[\Sigma Z, \Sigma Z]$ and the set $[\Omega \Sigma Z , \Omega \Sigma Z]_H$ of homotopy classes of H-maps $\Omega \Sigma Z \rightarrow \Omega \Sigma Z$ (see [22, p. 75]). Therefore, by taking the adjointness of the Samelson product, we obtain the result.

Step 2. Secondly, we suppose that the result holds for the $(k-1)$-fold iterated Whitehead product. Since $[\hat \rho_{i_{k-1}},\ldots , [\hat \rho_{i_1}, \hat \rho_{i_2}]\ldots]|_{Z_{4(i_1 + \cdots + i_{k-1})-1}} \simeq *$, and the map
$$
<\hat x_{i_{k-1}},\ldots , <\hat x_{i_1}, \hat x_{i_2}>\ldots>~~ : ~~S^{4(i_1 + \cdots + i_{k-1})} \longrightarrow \Omega \Sigma Z
$$
is rationally non-trivial due to Lemma 3.5 and adjointness again, by using Step 1 together with the maps, $\hat \rho_{i_k}$ and $\hat x_{i_{k}}$, we can construct an iterated commutator
$$
[\hat \rho_{i_k}, [\hat \rho_{i_{k-1}},\ldots , [\hat \rho_{i_1}, \hat \rho_{i_2}]\ldots ]] : Z \longrightarrow \Omega \Sigma Z
$$
such that, after taking the adjointness, the desiring formula holds.
\hfill\qquad$\square$
\enddemo

\bigskip

To complete the proof of the Theorem 3.1, we next observe that $\Sigma Z$ has the rational homotopy type of the bouquet of spheres; that is,
$$
\Sigma Z ~\simeq_{\Bbb Q}~ S^{5} \vee S^{9} \vee S^{13} \vee \cdots \vee  S^{4n+1} \vee \cdots,
$$
where $\simeq_{\Bbb Q}$ means a rational homotopy equivalence.
Thus, by using the salient Hilton's formula [12, Theorem A] for the homotopy groups of a wedge of spheres, we can find various kinds of rationally non-trivial indecomposable and decomposable generators on $\pi_{4n+1} (\Sigma Z) \otimes \Bbb Q$. For example, $\{\chi_1\}$ in dimension $5$, $\{\chi_2\}$ in dimension $9$, $\{\chi_3, [\chi_1 , \chi_2 ]\}$ in dimension $13$, $\{\chi_4,$ $[\chi_1 , \chi_3 ], [\chi_1, [\chi_1 , \chi_2 ]]\}$ in dimension $17$, $\{\chi_5 , [\chi_1 , \chi_4 ]$, $[\chi_1, [\chi_1 , \chi_3 ]], [\chi_1,$ $[\chi_1 [\chi_1 , \chi_2 ]]]$, $[\chi_2, \chi_3 ]$, $[\chi_2,$ $[\chi_1 , \chi_2 ]]\}$ in dimension $21$, and so on. Therefore, we can obtain the corresponding indecomposable and decomposable elements on $\pi_{4n+1}(\Sigma Z)/\text{torsion}$ for each $n = 1,2,3, \ldots$ (see the proof of Theorem 4.1).

Let $L$ and $L_{\leq n}$ denote the graded quasi-Lie algebras (Whitehead algebras) under Whitehead products on $\pi_{*}(\Sigma Z)/\text{torsion}$ and $\pi_{\leq 4n+1}(\Sigma Z)/\text{torsion}$ (corresponding to $\Cal L$ and $\Cal L_{\leq n}$), respectively. Then, for the triviality of the set, $SNT(\Sigma Z)$, it suffices to show that the map $\text{Aut} (\Sigma Z) \rightarrow \text{Aut}( L_{\leq n})$ has a finite cokernel for each $n$ since the map $\text{Aut}( \pi_{\leq 4n+1}(\Sigma Z)) \rightarrow \text{Aut}( L_{\leq n})$ has a finite kernel and it is an epimorphism.
In order to do this, we let $I_n L$ and $D_n L$ denote the indecomposable and decomposable components, respectively, of the group $\pi_{4n+1}(\Sigma Z)/\text{torsion}$. Then
\roster
\item $I_{n} L \cong \Bbb Z$ and thus $\text{Aut} (I_{n} L) \cong \Bbb Z_2$ for each $n$; and
\item the sequence $0 \rightarrow \text{Hom}(I_n L, D_n L) (\cong D_n L ) \rightarrow \text{Aut} (L_{\leq n}) \rightarrow
\text{Aut} (L_{< n}) \oplus \Bbb Z_2 \rightarrow 0$ is exact for each $n$ (see [18, p. 289]).
\endroster
Here the map $\text{Hom}(I_n L, D_n L) \rightarrow \text{Aut}(L_{\leq n })$ sends $[\rho_{i_k},[\rho_{i_{k-1}},\ldots,[\rho_{i_{1}}, \rho_{i_2}]\ldots]]_\sharp$ to $I + j \circ [\rho_{i_k},[\rho_{i_{k-1}},\ldots,[\rho_{i_{1}}, \rho_{i_2}]\ldots]]_{\sharp} \circ q$, and the map out of $\text{Aut}(L_{\leq n})$ is given by restriction and projection, where $q : L_{\leq n} \rightarrow I_n L$ is the projection and $j : D_n L \hookrightarrow L_{\leq n}$ is the inclusion. We note that the above exact sequence is still valid since we work on $\pi_{\leq 4n+1}(\Sigma Z)/\text{torsion}$. Moreover, we can see that
$\text{Aut} (\pi_{\leq 5} (\Sigma Z)/\text{torsion}) \cong {\Bbb Z}_2$ and $\text{Aut} (\pi_{\leq 9} (\Sigma Z)/\text{torsion}) \cong {\Bbb Z}_2 \oplus {\Bbb Z}_2$
and that $\text{Aut} (\pi_{\leq 4n+1} (\Sigma Z)/\text{torsion})$ is infinite for $n \geq 3$ (see Appendix). Thus the first induction step begins. Suppose that the map
$\text{Aut}(\Sigma Z) \rightarrow \text{Aut}(L_{< n})$ has a finite index.
We now discuss the self-homotopy equivalences which induce automorphisms on $\pi_{4n+1}(\Sigma Z)/\text{torsion}$. For each iterated Whitehead product $[x_{i_k},[x_{i_{k-1}},\ldots,[x_{i_{1}}, x_{i_2}]\ldots]] \in D_n L$ with $i_1 + i_2 + \cdots + i_k =n$ and $x_n \in I_n L$, by Lemma 3.6, we can construct a self-homotopy equivalence $I + [\rho_{i_k},[\rho_{i_{k-1}},\ldots,[\rho_{i_{1}}, \rho_{i_2}]\ldots]]$ in $\text{Aut}(\Sigma Z)$ which completely depends on the type of $[x_{i_k},[x_{i_{k-1}},\ldots,[x_{i_{1}}, x_{i_2}]\ldots]]$ such that the restriction
$$
(I + [\rho_{i_k},[\rho_{i_{k-1}},\ldots,[ \rho_{i_{1}}, \rho_{i_2}]\ldots]])_\sharp \vert_{L_{<n}}
$$
is the identity, and
$$
(I + [ \rho_{i_k},[\rho_{i_{k-1}},\ldots,[\rho_{i_{1}}, \rho_{i_2}]\ldots]])_\sharp (x_{n}) = x_{n} +  \alpha [x_{i_k},[x_{i_{k-1}},\ldots,[x_{i_{1}}, x_{i_2}]\ldots]],
$$
where $\alpha \neq 0$.
Moreover, there is only one indecomposable generator, up to sign, in $\pi_{4n+1}(\Sigma Z)/\text{torsion}$ for each $n$, while there are various kinds of decomposable generators in it (we denote these generators by the same notations for convenience); that is,
$$
\split
\pi_{4n+1} &(\Sigma Z)/\text{torsion} \\
&\cong \Bbb Z \{x_n\} \oplus
\displaystyle{\bigoplus_{[x_{i_k},[x_{i_{k-1}},\ldots,[x_{i_{1}}, x_{i_2}]\ldots]] \in D_n L }}
\Bbb Z\{ [x_{i_k},[x_{i_{k-1}},\ldots,[x_{i_{1}}, x_{i_2}]\ldots]] \}.
\endsplit
$$
Therefore, by using the induction hypothesis together with Theorem 3.2, we finally complete the proof of the
Theorem 3.1.  \hfill\qquad$\square$
\enddemo

\bigskip

In summary, the sets of all same $n$-type structures for the suspension of the classical infinite real, complex, and quaternionic projective spaces $\Sigma {\Bbb F}P^\infty ({\Bbb F} =  {\Bbb R}, {\Bbb C}, {\Bbb H})$ as well as the suspension of the Eilenberg-MacLane spaces of type $({\Bbb Z},n)$ have been constructed as follows:
Let $X ={\Bbb F}P^\infty$ or $K(\Bbb Z, n)$, where $\Bbb F =  \Bbb R, \Bbb C$ or $\Bbb H$. Then $SNT(\Sigma X)=*$ in the SNT-sense.

\bigskip

\heading {4. Appendix}
\endheading

\bigskip

By using the self-homotopy equivalences constructed in Remark 2.1, we obtain a concrete group structure of the automorphisms of the Whitehead algebras as follows:

\bigskip

\proclaim{Theorem 4.1} $\text{Aut} (\pi_{\leq n} (\Sigma Z )/\text{torsion})$ is finite for $n \leq 9$, and infinite for $n \geq 13$, and $\text{Aut}(\pi_{*} (\Sigma Z )/\text{torsion})$ becomes non-abelian. Moreover,
$\text{Aut}(\pi_{*} (\Sigma Z_{\Bbb Q} ))$ is both infinite and non-abelian.
\endproclaim

\demo{Proof}
We observe that there are infinitely many indecomposable generators and the iterated Whitehead products on $\pi_* (\Sigma Z)$ (resp. $\pi_* (\Sigma Z_{\Bbb Q})$) as follows:
\item{-} $x_1 , (\text{resp.}~ \chi_1 )$ in dimension $5$
\item{-} $x_2 , (\text{resp.}~ \chi_2 )$ in dimension $9$
\item{-} $x_3 , (\text{resp.}~ \chi_3 )$, and $[x_1 , x_2 ], (\text{resp.}~ [\chi_1 , \chi_2 ]),$ in dimension $13$
\item{-} $x_4 , (\text{resp.}~ \chi_4 )$, and $[x_1 , x_3 ]$, $[x_1, [x_1 , x_2 ]]$, $(\text{resp.}~ [\chi_1 , \chi_3 ], [\chi_1, [\chi_1 , \chi_2 ]])$ in dimension $17$
\item{-} $x_5 , (\text{resp.}~ \chi_5 )$, and $[x_1 , x_4 ]$, $[x_1, [x_1 , x_3 ]]$, $[x_1, [x_1 [x_1 , x_2 ]]]$,
$[x_2, x_3 ]$, $[x_2, [x_1 , x_2 ]]$ \newline $(\text{resp.}~ [\chi_1 , \chi_4 ]$, $[\chi_1, [\chi_1 , \chi_3 ]]$, $[\chi_1, [\chi_1 [\chi_1 , \chi_2 ]]]$, $[\chi_2, \chi_3 ]$, $[\chi_2, [\chi_1 , \chi_2 ]])$ in dimension $21$ \newline
$\ldots$. \newline
We thus have
\roster
\item $\text{Aut} (\pi_{\leq 5} (\Sigma Z)/\text{torsion}) \cong {\Bbb Z}_2$; and
\item $\text{Aut} (\pi_{\leq 9} (\Sigma Z)/\text{torsion}) \cong {\Bbb Z}_2 \oplus {\Bbb Z}_2$.
\endroster
If we take the self-map $\Psi  = I + [\rho_1 , \rho_2 ]$ of $\Sigma Z$, then the
induced homomorphism $\Psi_{\sharp}  = (I + [\rho_1 , \rho_2 ])_{\sharp}$ in $\text{Aut} (\pi_{\leq 13} (\Sigma Z )/\text{torsion})$ sends
$x_1 \mapsto x_1, x_2 \mapsto x_2$
and
$x_3 \mapsto x_3 + \alpha [x_1 , x_2]$, where $\alpha \neq 0$.
$\Psi_{\sharp}^2$ carries
$x_1 \mapsto x_1, x_2 \mapsto x_2$
and
$x_3 \mapsto x_3 + \alpha [x_1 , x_2] + \alpha [x_1 , x_2] + 0 \,\,\,( = x_3 + 2 \alpha  [x_1 , x_2])$.
Generally, by induction, $\Psi_{\sharp}^n$ sends $x_3 \mapsto x_3 + n \alpha [x_1 , x_2]$ for
each $n$. Thus $\Psi_{\sharp}$ has an infinite order in $\text{Aut} (\pi_{\leq n} (\Sigma Z )/\text{torsion})$ for $n \geq 13$.

In order to show that $\text{Aut}(\pi_{*} (\Sigma Z )/\text{torsion})$ is non-abelian, it suffices to show that
$\text{Aut} (\pi_{\leq 4m+5} (\Sigma Z )/\text{torsion}), m \geq 3$ is non-abelian: Let $f$ and $g$ be self-maps of $\Sigma Z$ inducing
$$
f_\sharp = {\cases (I + [\rho_1 , \rho_{m-1} ])_\sharp
  &\text{for dim $\leq 4m+1$}, \\
      I_\sharp &\text{otherwise} \endcases}
$$
and
$$
g_\sharp = {\cases (I + [\rho_1 , \rho_m ])_\sharp &\text{for dim $\leq 4m+5$}, \\
      I_\sharp &\text{otherwise} \endcases}
$$
between homotopy groups, where `dim' means the homotopy dimension. Then $f_\sharp$ sends
$$
\split
&x_1 \longmapsto x_1 , \hskip 1.5cm x_2 \longmapsto x_2 , \hskip 0.5cm \ldots , \hskip 0.5cm  x_{m-1} \longmapsto x_{m-1},\\
&x_m \longmapsto x_m + \alpha_1 [x_1 , x_{m-1}] ,\hskip 1.0cm  x_{m+1} \longmapsto x_{m+1}, \hskip 0.5cm  \ldots,\\
\endsplit
$$
where $\alpha_1 \neq 0$.
Moreover, $g_\sharp$ sends
$$
\split
&x_1 \longmapsto x_1 , \hskip 0.7cm  x_2 \longmapsto x_2 , \hskip 0.5cm \ldots , \hskip 0.5cm x_{m-1} \longmapsto x_{m-1},\\
&x_m \longmapsto x_m ,  \hskip 0.4cm x_{m+1} \longmapsto  x_{m+1} + \alpha_2 [x_1 , x_m ], \\
&x_{m+2} \longmapsto x_{m+2}, \hskip 0.5cm \ldots ,
\endsplit
$$
where $\alpha_2 \neq 0$.
We note that the above generators $x_{m+2} (m \geq 3)$ might be killed off in
$$
\text{Aut}(\pi_{\leq 4m+5} (\Sigma Z )/\text{torsion})
= \text{Aut}(\pi_{*} (\Sigma Z )^{(4m+5)}/\text{torsion}),
$$
where $(\Sigma Z )^{(4m+5)}$ means the Postnikov approximation of $\Sigma Z$. We thus have
$$
\split f_\sharp \circ g_\sharp (x_{m+1}) &= f_\sharp (x_{m+1} + \alpha_2 [x_1 , x_m ])\\
                        &= x_{m+1} + \alpha_2 [x_1 , x_m +\alpha_1 [x_1 ,x_{m-1} ]] \\
                        &= x_{m+1} + \alpha_2 [x_1 , x_m ] + \alpha_1 \alpha_2 [x_1 , [x_1 , x_{m-1}]]\\
\endsplit
$$
and
$$
g_\sharp \circ f_\sharp (x_{m+1} ) = g_\sharp (x_{m+1} ) = x_{m+1} + \alpha_2 [x_1 , x_m ].
$$
We note that both $[x_1 , x_m ]$ and $[x_1 , [x_1 , x_{m-1}]]$ are basic Whitehead products. Thus by the Hilton's theorem [12] they are nonzero and cannot be equal to each other. Therefore $f_\sharp \circ g_\sharp \neq g_\sharp \circ f_\sharp$.

Similarly, we can show that $\text{Aut}(\pi_{*} (\Sigma Z_{\Bbb Q} ))$ is both infinite and non-abelian by considering the action of self-homotopy equivalences of the form
$$
I + [\rho_{i_k},[\rho_{i_{k-1}},\ldots,[\rho_{i_{1}}, \rho_{i_2}]\ldots]]
$$
on the set of indecomposable generators $\{\chi_i\}_{i \in \Bbb Z_+}$ of rational homotopy groups. Indeed,
\roster
\item $\text{Aut}(\pi_{\leq 5} (\Sigma Z_{\Bbb Q})) \cong (\Bbb Q - \{0\})$,
\item $\text{Aut}(\pi_{\leq 9} (\Sigma Z_{\Bbb Q})) \cong (\Bbb Q - \{0\}) \oplus (\Bbb Q - \{0\})$,
\endroster
and the rest of the proof goes to the same method in the case of $\text{Aut}(\pi_{*} (\Sigma Z )/\text{torsion})$, a group of automorphisms of integral homotopy modulo torsion. \hfill\qquad $\square$
\enddemo

\bigskip

There is a theorem, originally due to Quillen, which states that any connected graded Lie algebra over $\Bbb Q$ is realized as $\pi_* (\Omega X)\otimes \Bbb Q$ for some $1$-connected space $X$. In our case, $\pi_* (\Sigma Z_{\Bbb Q})$ is free graded quasi-Lie algebra on generators $\chi_1, \chi_2,\ldots, \chi_n, \ldots$ in rational homotopy theory.

What will happen in the case of the infinite complex projective space $\Bbb CP^\infty$?
We note that it is an H-space, while $Z = \Bbb HP^\infty$ is a rational H-space as just mentioned before. Similarly, just as in the case of $\Bbb HP^\infty$, we have the following corollary.

\bigskip

\proclaim{Corollary 4.2} $\text{Aut} (\pi_{\leq n} (\Sigma {\Bbb C}P^\infty )/\text{torsion})$ is finite for $n \leq 5$ and it is an infinite non-abelian group for $n \geq 7$. Furthermore, $\text{Aut}(\pi_{*}( \Bbb CP^\infty_{\Bbb Q} ))$ are both infinite and non-abelian groups.
\endproclaim

\demo{Proof} The classical Hopf-Thom theorem and the Serre spectral sequence of a fibration assert that
$$
\Sigma \Bbb CP^\infty \simeq_{\Bbb Q} S^3 \vee S^5 \vee \cdots \vee S^{2i+1} \vee \cdots,
$$
where $\simeq_{\Bbb Q}$ means the rational homotopy equivalence.
Thus we get the corresponding indecomposable generators and iterated Whitehead products on homotopy groups, say,
$$
\xi_3 , \xi_5 , \xi_{7}, [\xi_3 , \xi_5 ], \xi_{9}, [\xi_3 , \xi_{7}], [\xi_3 , [\xi_3 , \xi_{5}]], \xi_{11} , \ldots,
$$
where $\xi_n$ has the dimension $n$ in homotopy groups of $\Sigma \Bbb CP^\infty$. The rest of the proof goes to the same way as in the one of Theorem 4.1. \hfill\qquad $\square$
\enddemo

\bigskip

Department of Mathematics, and Institute of Pure and Applied Mathematics, Chonbuk National University, Jeonju, Jeonbuk 561-756, Republic of Korea

Tel.: +82-63-270-3367, Fax: +82-63-270-3363

E-mail address: dwlee\@jbnu.ac.kr

\enddocument